\documentclass[10pt, conference]{IEEEtran}
\IEEEoverridecommandlockouts
\usepackage{soul}
\usepackage{xcolor}
\usepackage[ruled,vlined,linesnumbered]{algorithm2e}
\usepackage{color, colortbl}
\definecolor{LRed}{rgb}{1,.8,.8}
\definecolor{MRed}{rgb}{0.8,.8,1}
\definecolor{HRed}{rgb}{1,.2,.2}
\usepackage{cite}
\usepackage{amsmath,amssymb,amsfonts}
\usepackage{algorithmic}
\usepackage{textcomp}
\usepackage{xcolor}
\usepackage{comment}
\usepackage{indentfirst}

\usepackage{multirow}
\usepackage{booktabs}
\usepackage{multirow}
\usepackage{amsthm}
\usepackage{float}
\usepackage{url}
\usepackage[pdftex]{graphicx}
\usepackage{subcaption}

\newcommand{\brac}[1]{\left\{ #1 \right\}}

\newcommand{\sZ}{\mathbb{Z}_{+}}
\newcommand{\sR}{\mathbb{R}}

\DeclareMathOperator*{\argmax}{argmax}
\DeclareMathOperator*{\argmin}{argmin}

\IEEEpubid{\makebox[\columnwidth]{978-1-6654-4177-
3/21/\$31.00~\copyright{}2021 IEEE \hfill}
\hspace{\columnsep}\makebox[\columnwidth]{ }}

\begin{document}
\title{Safe Sequential Optimization for \\Switching Environments}
\author{Durgesh Kalwar and Vineeth B. S. \thanks{The first author is with the Department of Mathematics, Indian Institute of Space Science and Technology, Thiruvananthapuram. The second author is with the Department of Avionics, Indian Institute of Space Science and Technology, Thiruvananthapuram. Emails: dghkalwar007@gmail.com, vineethbs@gmail.com}}

\maketitle

\begin{abstract}
We consider the problem of designing a sequential decision making agent to maximize an unknown time-varying function which switches with time.
At each step, the agent receives an observation of the function's value at a point decided by the agent.
The observation could be corrupted by noise.
The agent is also constrained to take safe decisions with high probability, i.e., the chosen points should have a function value greater than a threshold. 
For this switching environment, we propose a policy called Adaptive-SafeOpt and evaluate its performance via simulations.
The policy incorporates Bayesian optimization and change point detection for the safe sequential optimization problem.
We observe that a major challenge in adapting to the switching change is to identify safe decisions when the change point is detected and prevent attraction to local optima.
\end{abstract}

\section{Introduction}
Safe optimization of unknown functions arises in many real-world scenarios such as robotic systems, unmanned exploratory vehicles, and autonomous cars.
For example, Krause et al. \cite{krause2017safe} considered the problem of an autonomous rover exploring the surface of Mars.
The height and gradient of the surface is unknown to the rover.
Since the rover has physical limitations with respect to the gradients it can move over, it has to safely explore the surface while visiting points that maximize scientific insight.
The problem of safe movement of the rover while optimizing scientific insight is an example of safe optimization of an unknown function.
Although the function is unknown in such settings, an optimizing agent can interact with the function and obtain a noisy observation of the function at a chosen point. 

Motivated by such applications, we consider the problem of designing an agent to safely find the maximum of an unknown time-varying function from noisy observations.
Each observation is a noisy evaluation of the function's value at a point in the domain chosen by the agent.
The points chosen by the agent are also required to meet a safety criterion.
Prior work \cite{sui2015safe}, \cite{wachi2018safe} has used the framework of Bayesian optimization in order to propose sequential decision making agents for the above safe optimization problem with time invariant functions.
Using standard reinforcement learning terminology, we model the optimizing agent as interacting with an environment.
The agent interacts with the environment by choosing points in the function domain.
The environment provides the agent with a noisy observation of the function value at the chosen point.

In this paper, we consider the important extension of the above problem to switching environments.
A switching environment is one in which the unknown function is time varying and exhibits a discontinuous \emph{switching} to another unknown function at a change-point epoch.
The non-stationarity of the environments that agents have to contend with is an important challenge for real-world problems \cite{49999}.
We propose a heuristic policy (which is an extension of the algorithm in \cite{sui2015safe} with change detection for the unknown functions) and evaluate the performance of the policy using simulations.
We observe that an important challenge in this problem is re-initializing an estimate of the safe set once the change has been detected and proposing a solution.

Bayesian optimization is used for addressing the problem of safe exploration and optimization in which the unknown objective function and/or safety constraints are modelled using Gaussian processes.
The quantification of uncertainty, which is obtained for free with the Bayesian framework, is used to decide whether an action is safe.
For bandit setting, Sui et al. \cite{sui2015safe} proposed Safe-Opt algorithm (Safe exploration for Optimization).
An apriori unknown safety function is modelled using Gaussian Processes (GP) and it's confidence interval is used to decide whether a sequence of decisions taken during exploration is safe or not.
If at a decision time the value of the safety function is more than a threshold then it is safe. 
In their setting the safety and objective function are identical.
The proposed Safe-Opt algorithm trades off between maximizing the size of the safe set of decisions starting from an initial safe seed and finding the optimal reachable decision in that safe set. 
We note that the authors considered that the unknown function is stationary with time.
In contrast to their work we consider a non-stationary scenario and propose a change point detection based extension to Safe-Opt.
We note that in addition to the tradeoff between maximising the size of the safe-set and optimal reachable decision we also have another tradeoff in the exploration required to detect the change-point.

For  Markov Decision Process (MDP) setting, Krause et al. \cite{krause2017safe} proposed a safe exploration algorithm called Safe-MDP. 
In their work they assumed that the transition model is known, but the safety function is unknown.
The safety function, which is assumed to have similar values in similar states, is then modelled using GP.
In their work they considered exploration in stationary MDP setting.
Wachi et al. \cite{wachi2018safe} proposed a safe exploration with an  optimization algorithm for finite deterministic MDP and provides theoretical guarantees on the satisfaction of the safety constraint.
However, the acquired policy is not necessarily near-optimal in terms of the cumulative reward. 
Wachi and Sui \cite{wachi2020safe} proposed a safe RL algorithm for finite deterministic MDP that guarantees a near-optimal cumulative reward while guaranteeing the satisfaction of the safety constraint as well. 
In their work they also assumed the transition model is known and both reward and safety function are modelled using GP. 
Wachi et. al. \cite{wachi2018safetime} extended Safe-MDP to the case of time-variant safety functions, which are assumed to be Lipschitz continuous with respect to time.

We note that optimization of unknown time varying (switching) functions without safety constraints has been addressed by many authors.
Mellor and Shapiro  \cite{mellor2013thompson} had proposed a Bayesian online change point detection based method for switching bandits.
A similar approach was also used by \cite{alami2017memory}.
Recently, Ghatak \cite{ghatak2020change} had proposed a change detection based Thompson sampling framework for non-stationary bandits.
Padakandla et al. \cite{padakandla2020survey} provides a survey of reinforcement learning algorithms for dynamically varying environments.
We note that our work incorporates the notion of safety in addition to the non-stationarity considered in the above papers.

We note that there are multiple approaches to ensuring \emph{safety} for agents which include the one summarized above.
Garcıa and Fernández \cite{garcia2015comprehensive} provides a succinct survey on safe reinforcement learning.
In the risk sensitive approach, the long-term reward maximization is modified to include risk measures, such as variance or higher moments of reward. However, these approaches only minimize risk and do not treat safety as a hard constraint. 
Alternatively, the optimization criterion is transformed to include the probability of visiting error states (e.g. Geibel and Wysotzki \cite{geibel2005risk}).
In other work on safe reinforcement learning, Moldovan and Abbeel  \cite{moldovan2012safe} consider the problem of safe exploration in MDPs. 
They ensure safety by restricting policies to be ergodic with high probability, i.e., able to recover from any state visited. 
This is computationally demanding even for small state spaces and doesn't provide convergence guarantees. 
Biyik et al. \cite{biyik2019efficient}, consider the problem of safe exploration in deterministic MDPs with unknown transition models. 
They considered safety criterion similar to that in \cite{moldovan2012safe}.
Roderick et al. \cite{roderick2020provably}, consider the problem of safe exploration in PAC (probably approximately correct)-MDP with unknown, stochastic dynamics. 

\noindent\textbf{Outline and Contributions}
We define the system model and our problem statement in Section \ref{subsec:System Model and Problem Statement}.
The Safe-Opt algorithm and related notation is then presented in Section \ref{sec:safeopt}.
The main contribution in this paper is the proposal of an algorithm Adaptive-SafeOpt for safe optimization of an apriori unknown function.
This algorithm is presented in Section \ref{sec:adaptivesafeopt}.
In Section \ref{sec:genie} we consider ``genie'' algorithms which serve as baselines to compare and understand the performance of Adaptive-SafeOpt.
Simulation results and discussions are presented in Section \ref{sec:performanceanalysis}.

\section{System Model and Problem Statement}
\label{subsec:System Model and Problem Statement}
We consider a discrete time model for the safe optimization problem.
The time epochs at which the function is evaluated is denoted as $t \in \sZ$.
Let $f(x, t)$ be a scalar valued function of $x$ and $t$, where $x \in \mathcal{X} \subseteq \sR$.
We assume that $\mathcal{X}$ is compact and $f(x,t)$ is Lipschitz continuous with respect to $x$.
The Lipschitz constant is assumed to be $L$.
The objective of the safe optimization problem is to find the maximum of $f(x,t)$ subject to safety constraints.
In this paper we consider switching environments, i.e., $f(x,t)$ is assumed to be $f_{1}(x)$ until an arbitrary change time $t_{c}$ and then $f(x, t) = f_{2}(x)$ for $t \geq t_{c}$.
We consider the problem with only one change\footnote{We note that the algorithms proposed in this paper can be extended to the case of multiple changes without any change.}.
We also assume that $\forall x$, $|f_{1}(x) - f_{2}(x)| \leq B$, where $B$ is positive.

We note that since the function is unknown, we optimize the function by observing the value of the function at points $x_{t} \in \mathcal{X}$ which are chosen for every time $t$.
The observed value at time $t$ is denoted as $y_{t}$.
The safety constraint that we consider in this paper is that the fraction of time $f(x_{t},t) \geq h$ is greater than or equal to $1 - \delta$ where $\delta \in (0, 1)$.
Here we note that $h \in \sR$ is chosen to be less than $\max_{x} f(x, t), \forall t$.
So ideally our problem is to obtain a sequence of points $x^{*}_{t}$ such that $\forall t, f(x^{*}_{t}, t) = \max_{x \in \mathcal{X}} f(x, t)$, and $f(x^{*}_{t}) \geq h$.
Since the functions are unknown we note that obtaining such as sequence of points would be a tough task.
We therefore first define the following metrics and formulate a simpler problem.

A policy $\pi$ is defined to be a sequence of $x_{t}$ chosen by the optimizing agent.
The sequence $x_{t}$ could be chosen as a function of the history of choices, i.e., $(x_{0},x_{1},\dots,x_{t - 1})$ as well as the observations $(y_{0}, y_{1}, \dots, y_{t-1})$.
We denote by $\Delta_{\pi}(t)$ the gap between the maximum value of the function and the observed value $y_{t}$ at $x_{t}$
\begin{eqnarray*}
	\Delta_{\pi}(t) = \max_{x \in \mathcal{X}} f(x, t) - y_{t}
\end{eqnarray*}
We also define the normalized cumulative regret over a horizon $T$ as
\begin{eqnarray*}
	R_{\pi}(T) = \frac{\sum_{t} \Delta_{\pi}(t)}{T}
\end{eqnarray*}
The cumulative unsafe evaluations over a horizon $T$ is defined as
\begin{eqnarray*}
    U_{\pi}(T) = {\sum_{t} \mathbb{I}\brac{f(x_{t}, t) < h}},
\end{eqnarray*}
where $\mathbb{I}$ is the indicator function.
Our objective in this paper is to find a policy $\pi$ such that the regret is minimized for all $T$ subject to a safety constraint.
\begin{eqnarray}
	\min_{\pi} R_{\pi}(T) \text{ such that } \frac{U_{\pi}(T)}{T} \leq \delta.
	\label{eq:problem}
\end{eqnarray}
We define the true safe set $\mathcal{S}^{*}_{t}$ as $\mathcal{S}^{*}_{t} = \brac{x: f(x, t) \geq h}$.
We assume that at $t = 0$, an element of $\mathcal{S}^{*}_{0}$ called safe-seed is known to the algorithm\footnote{For comparison purposes, we introduce genie policies. For these genie policies we note that at $t = 0$, one point each from the disjoint intervals that makes up $\mathcal{S}^{*}_{0}$ is given as input. This is called the safe-seed set.}.

Ideally, we would want a policy $\pi$ that achieves the above minimum for any choice of $f_{1}$ and $f_{2}$.
However, this may not be possible \cite[Chapter 2]{10.1145/3440959.3440965}.
In this paper, we evaluate the performance of a policy by considering the average of the above metrics over randomly chosen pairs $(f_{1}, f_{2})$.

\section{Background: Safe-Opt Algorithm}
\label{sec:safeopt}
In this section, we briefly review the Safe-Opt algorithm proposed by Sui et. al. \cite{sui2015safe} and introduce some essential notation since our algorithm is an extension of Safe-Opt to switching environments.
For the time-invariant environment in \cite{sui2015safe}, the maximization of the unknown function $f(x)$ is done by estimating $f(x)$ using Gaussian process (GP) regression.
In GP regression, the unknown function is assumed to be modelled by a sample function from a GP prior \cite{williams2006gaussian}.
The GP prior is completely characterized by its mean function $\mu(x)$ (without loss of generality $\mu(x) = 0$) and covariance function $k(x, x')$ where $x, x' \in \mathcal{X}$.
At every time $t$, the Safe-Opt policy chooses a point $x_{t}$ and receives an observation $y_{t} = f(x_{t}) + n_{t}$ where $n_{t}$ is an independent sample from Gaussian noise with mean $0$ and variance $\sigma^{2}$.
Based on $y_{t}$ a posterior distribution for the unknown function can be derived.
This posterior distribution is again Gaussian and characterized completely by a mean function $\mu_{t}(x)$ and covariance function $k_{t}(x, x')$.
In order to satisfy the safety constraints, Safe-Opt computes upper and lower confidence bounds on the function using this posterior.
The upper $u_{t}(x)$ and lower $l_{t}(x)$ confidence bounds are defined as
\begin{eqnarray*}
	u_{t}(x) & = & \mu_{t}(x) + \beta_{t} \sigma_{t}(x), \\
	l_{t}(x) & = & \mu_{t}(x) - \beta_{t} \sigma_{t}(x).
\end{eqnarray*}
By an appropriate choice of $\beta_{t}$ and by choosing $x_{t + 1}$ such that $l_{t}(x_{t + 1}) \geq h$, Safe-Opt is able to satisfy the safety constraint with high probability.
We denote by $Q_{t}(x)$ the interval $[l_{t}(x), u_{t}(x)]$ as a function of $x \in \mathcal{X}$.
We also denote the length of the confidence interval as $w_t(x):=u_t(x) - l_t(x)$.
We note that on the basis of the confidence bounds an estimate $S_{t} \subseteq \mathcal{X}$ of the safe set can be maintained which is defined as
\begin{equation}
	S_t = \brac{x \in \mathcal{X} | l_{t}(x) \geq h}
	\label{eqn: safe set}
\end{equation}
We note that at $t = 0$ we are given a safe seed $\in \mathcal{S}^{*}_{0}$ (in the context of Safe-Opt $\mathcal{S}^{*}_{t}$ is time-invariant).
Since the safe seed may not achieve the maximum of $f(x)$ we need to explore safely.
Safe-Opt maintains a set $G_t \subseteq S_t$ of candidate decisions that, upon potentially repeated selection, have a chance to expand $S_t$. The set $G_t$ is defined as 
\begin{equation}
G_t = \{ x \in S_t | \psi_t (x) > 0\}
\end{equation}
where
\begin{equation*}
\psi_t(x)=\vert \{  x' \in \mathcal{X} \setminus S_t| u_t(x)-Ld(x, x' ) \geq h\}\vert.
\end{equation*}
We note that SafeOpt assumes Lipschitz continuity for the function $f(x)$ with Lipschitz constant $L$ over $x \in \mathcal{X}$.
We also note that in order to find the maxima, we need to consider candidate points which are chosen from a set  $M_t \subseteq S_t$ of decisions that are potential maximizers of $f$.
\begin{equation}
M_t = \{x \in S_t | u_t (x) \geq \max_{x' \in S_t} l_t(x')\}
\end{equation}
Safe-Opt policy then chooses points $x_{t}$ according to
\begin{equation}
	x_t = \argmax_{x \in G_t \cup M_t }w_t(x).
\end{equation}

\section{Adaptive SafeOpt}
\label{sec:adaptivesafeopt}
In this section, we propose a heuristic policy (Adaptive SafeOpt) that extends Safe-Opt \cite{sui2015safe} to adapt to the switches in $f(x,t)$.
We note that an intuitive approach to adapting to the change in the function $f(.)$ safely is to detect whether a change has happened and then restart the Safe-Opt algorithm with a new safe seed.
The challenges here are therefore to quickly detect the change as well as to find a safe seed for restarting Safe-Opt.
In contrast to Safe-Opt, Adaptive SafeOpt balances three objectives: the desire to expand the safe region, the need to obtain $x_{t}$ which achieves the maxima, and the need to detect the change-point.

We note that the following is a candidate rule which can be used to detect a change.
At each time step $t$ we observe a noisy observation of function $f$, $y_t=f(t,x_t)+n_t$, from which we update the GP model of function, where $x_t$ is sampled according to the above sampling criteria. 
To detect the change-point, at every time step we check the condition that the observed $y_t$ is within the current confidence interval $Q_{t}$ or not.
If $y_t \in Q_{t}(x_{t})$ then the algorithm decides that the function has not changed.
If $y_t \not \in Q_{t}(x_{t})$ then Adaptive-SafeOpt declares that the change-point has detected and the function has changed.
In order to balance between the need to detect a change as well as maximize the function safely, we use an $\epsilon$-greedy approach for Adaptive SafeOpt.
At every time $t$ we choose
\begin{equation}
	x_t=
	\begin{cases}
		\argmin_{x \in S_t  }w_t(x) & \text{with $\epsilon$ probability }  \\
		\argmax_{x \in G_t \cup M_t }w_t(x) & \text{with $1-\epsilon$  probability }  \\
	\end{cases}
\end{equation} 
Suppose a change has been detected, then we also need to estimate a new safe set $S_{t}$.
If the $y_{t}$ at the declared change time is safe, then the new safe seed is $x_{t}$ itself.
On the other hand if $y_{t} < h$, then we initialize a safe-set estimate defined as
\begin{equation}
S_t = \{x \in \mathcal{X} | l_{t-1}(x) - B \geq h\}
\label{eqn: safe set at change}
\end{equation}
Here we make use of the assumption that $|f_{1}(x) - f_{2}(x)| \leq B, \forall x$.
It may turn out that $S_{t} = \emptyset$ or not.
If $S_{t} \not = \emptyset$ then we have a safe set and we continue with Safe-Opt as before.
However, if $S_{t} = \emptyset$ according to the above rule then we pick a $x_{t + 1}$ from $\argmax l_{t}(x)$. 
We note that $x_{t}$ has been used to update the GP, although it is unsafe.

\begin{algorithm}[!ht]
	\SetAlgoLined
	\scriptsize
	\caption{Adaptive-SafeOpt}
	\label{algo:adaptive window}
	\KwIn{Function domain $\mathcal{X}$, GP prior $(\mu, k)$, signal variance parameter $\sigma_0$, seed set $S_0$, safety threshold $h$, $window\_min$, $window\_max$, $delaychangedetection\_flag$ = True, $changepoint\_flag$ =False, $changedetection\_delay$, $counter = 0$, $changepoint\_index = 1$, $B$, $\epsilon$.}
	Initialize GP with safe seed points $S_0$ and compute $Q_{0}$\\
	$X=\{x | x \in S_0\}$, $Y=\{f(x)|x \in S_0\} $ \\
	\For {$t = 1, . . $}
	{ 
		\eIf{changepoint\_flag = false}
			{
				$S_t \leftarrow \{x \in \mathcal{X} | l_t(x) \geq h\}$\\
				$M_t \leftarrow \{x \in S_t | u_t (x) \geq \max_{x' \in S_t} l_t(x')\}$\\
				$G_t \leftarrow \{ x \in S_t | \psi_t (x) > 0\}$\\
			}
		{
			$S_t \leftarrow \{x \in \mathcal{X} | l_{t-1}(x)-B \geq h\}$\\
			$M_t \leftarrow \{x \in S_t | u_{t-1}(x) \geq \max_{x' \in S_t} l_{t-1}(x')\}$\\
			$G_t \leftarrow \{ x \in S_t | \psi_{t-1} (x) > 0\}$\\
			$changepoint\_flag = False$
		}
		\eIf{delaychangedetection\_flag = True}
			{
				$$x_t \leftarrow 
				\begin{cases}
				\argmax_{x \in G_t \cup M_t }(w_t(x)) & \text{if $S_t \neq \emptyset $}\\
				\argmax_{x \in \mathcal{X}} (l_t(x)) & \text{if $S_t = \emptyset$}
				\end{cases}$$\\
				$y_t \leftarrow f(x_t ) + n_t$\\
				$window=window+window\_increment$\\
				$counter=counter+1$\\
				\If{$window > window\_max$}
				    {
				    $window = window\_max$
				    }
				\If{counter = changedetection\_delay}
				{
				    $counter=0$ \\
				    $delaychangedetection\_flag=False$.
				}
			}
		{
			$$x_t \leftarrow
			\begin{cases}
			\argmin_{x \in S_t  }w_t(x) & \text{with $\epsilon$ probability }  \\
			\argmax_{x \in G_t \cup M_t }w_t(x) & \text{with $1-\epsilon$  probability }  \\
			\argmax_{x \in \mathcal{X}} (l_t(x)) & \text{if $S_t = \emptyset$}
			\end{cases}
			$$\\
			$y_t \leftarrow f(x_t ) + n_t$\\
			\eIf{$y_t < l_t(x_t)$ or $y_t > u_t(x_t)$}
			{   
				$window = window\_min$\\
				$changepoint\_index = t$\\
				$delaychangedetection\_flag = True$\\
				$changepoint\_flag = True$\\
			}
			{
				$window=window+window\_increment$\\
				\If {$ window > window\_max$ }
				 {
				 $ window = window\_max $
				 }
			}
		}
		start=t-window\\
		\If {start $<$ changepoint\_index} 
		{ $start = changepoint\_index$ }
		Update GP using $(x_{start}, \dots, x_{t})$ and $(y_{start}, \dots, y_{t})$.
		Compute $Q_t(x), \forall x \in S_t$ \\
	}
\end{algorithm}

The complete algorithm is given in Algorithm \ref{algo:adaptive window}.
The notation used in Algorithm \ref{algo:adaptive window} is defined in Section \ref{sec:safeopt}.
A few practically motivated modifications are also used in Algorithm \ref{algo:adaptive window}.
Suppose we have prior information about the inter-change duration, e.g., we know that the inter-change duration is at least some number of slots.
Then, the $\epsilon$-greedy policy need not be used immediately after a change-point.
We incorporate this by not using the above $\epsilon$-greedy policy until a counter expires.
In order to reduce data storage, we also introduce a data window. 
The data window size is incremented by one until a maximum window size ($window\_max$) is reached.

We note that an intuitive method to handle a time-variant environment is to consider data only in the immediate past. 
In order to evaluate how the Adaptive-SafeOpt policy compares with such a policy we also consider a FixedWindow-SafeOpt policy defined as follows.
The FixedWindow-SafeOpt policy has a parameter \emph{window}.
For FixedWindow-SafeOpt, the GP model for $f(x,t)$ is updated at every time $t$ using $(x_{t - window + 1}, \dots, x_t)$ and $(y_{t - window + 1},y_t)$.
Then the sets $Q_{t}(x), G_{t}$, and $M_{t}$ are computed and $x_{t + 1}$ is chosen as in SafeOpt (see Section \ref{sec:safeopt}).

\section{Algorithms for comparison}
\label{sec:genie}
In this section we discuss ``Genie'' algorithms which have access to \emph{extra} or \emph{side} information.
Genie policies are not practically implementable since they assume the availability of such information, but are used as baselines for comparing the performance of implementable policies such as Adaptive-SafeOpt.

\noindent{\textbf{Genie-CP-SS}}: This is a policy that has knowledge of the time $t_{c}$ at which change point happens as well as the true safe seed set for $f_{2}$  after switching.
We note that a function ($f_{1}$ or $f_{2}$) may have multiple disjoint intervals in the true safe set. 
We assume that a single point from each of these disjoint intervals is given as part of the safe seed set to Genie-CP-SS.
Then, for $t < t_{c}$ Genie-CP-SS uses Safe-Opt which is initialized with the safe seed, and for $t \geq t_{c}$ Safe-Opt can be re-initialized with the new safe seed and used.
Thus, the policy chooses $x_t = \argmax_{x \in G_t \cup M_t }(w_t(x))$.
We note that since $t_{c}$ as well as the safe-seed set is known, Genie-CP-SS should achieve the minimum possible value of regret with the minimum number of unsafe evaluations and provides an useful baseline for comparing with Adaptive-SafeOpt.

\noindent{\textbf{Genie-CP}}: This policy has side information only about the change point and not about the safe seed when a change happens.
At $t_{c}$ if $y_{t_{c}} \geq h$, then we re-initialize $S_{t_{c}} = x_{t_{c}}$
 and then for $t > t_{c}$, the policy chooses $x_{t}$ according to $x_t = \argmax_{x \in G_t \cup M_t }(w_t(x))$.
Otherwise, we choose $x_{t}$ as $\argmax_{x \in \mathcal{X}} l_t(x)$.
The performance of Genie-CP would indicate the loss in performance due to the non-knowledge of safe seed set.

\noindent{\textbf{Genie-SS}}: This policy has side information of the safe seeds. However, it does not know the change point and uses a change point detection scheme as follows (this is similar to used by Adaptive Safe-Opt).
At every time $t$, if the algorithm is allowed to do change detection (see discussion about incorporating prior information about change point times for Adaptive-SafeOpt), and if the current observation $y_{t} \not \in Q_{t}(x_{t})$ Genie-SS declares that a change has happened.
Once a change is declared to have happened, the genie is given one safe seed each from each of the disjoint intervals which makes up the true safe set $S^{*}_{t}$.
Similar to Genie-CP, the performance of Genie-SS would indicate the loss in performance due to non-knowledge of the change point.

\noindent{\textbf{GP-UCB-CP}}: Srinivas et al. \cite{Srinivas10gaussianprocess} had proposed GP-UCB algorithm which does not consider the safety constraint. Here we consider GP-UCB endowed with side information about when the change occurs so as to compare the regret with our policy. We note that GP-UCB-CP uses $x_t = \argmax_{x \in \mathcal{X}} u_t(x)$ and re-initializes the algorithm at the change point $t_{c}$.
We note that GP-UCB-CP would inform us about the global maxima without any regard to the safety constraint.
We expect that Genie-CP-SS and GP-UCB-CP would perform similarly as the global maxima is also safe; however we note that the exploration methodology for both of these policies are different.

\section{Simulations and Performance Analysis}
\label{sec:performanceanalysis}
\begin{figure}[ht!]
    \centering
    \includegraphics[width=90mm, height=45mm]{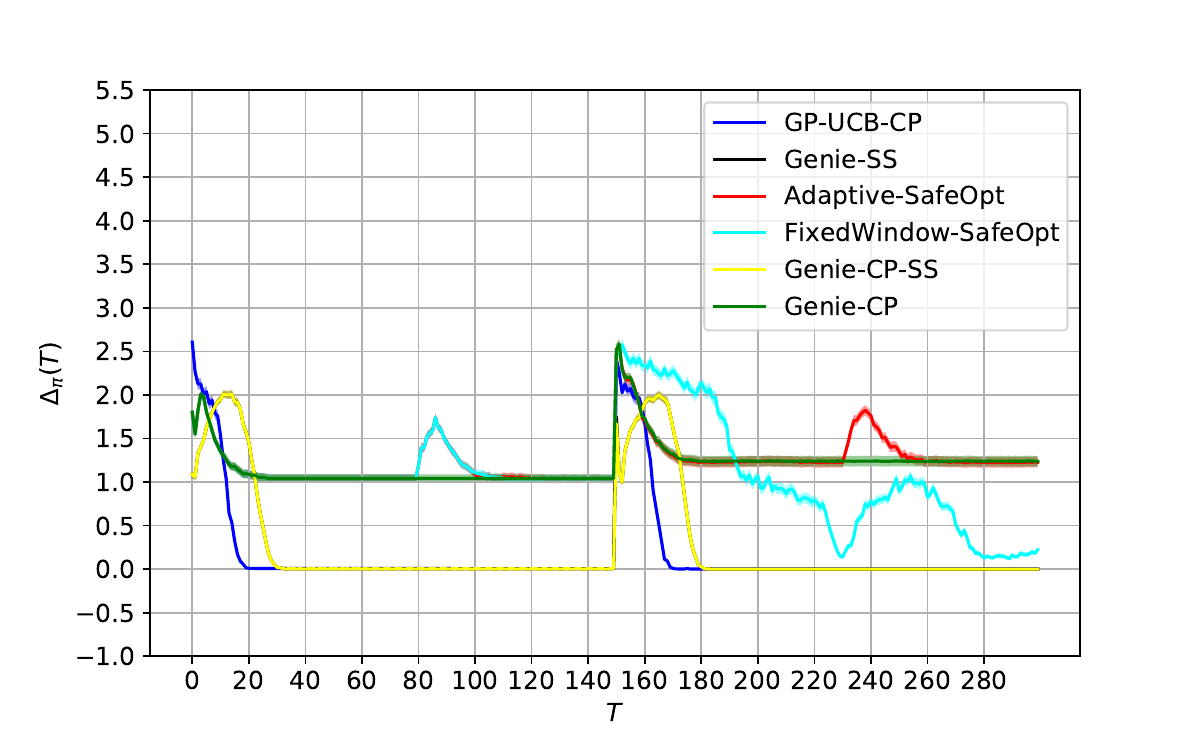}
    \caption{Comparison of $\Delta_{\pi}(T)$ as a function of $T$ for different algorithms. The change point $t_{c} = 150$. In this illustration, we assume that there is no observation noise.}
    \label{fig:pseudo regert with outlier}
\end{figure}

For comparing the performance of the algorithms proposed above, we consider one-dimensional functions $f_{1}(x)$ and $f_{2}(x)$ which are sampled from a GP prior.
The safety threshold $h$ is assumed to be $0$ without loss of generality.
The mean function $\mu(x)$ is assumed to be $0$ and the covariance function is specified by a radial-basis function kernel (parameterized by variance of $2$ and length scale of $1$).
When sampling $f_{2}$ we restrict to those samples such that $\forall x, |f_{1}(x) - f_{2}(x)| \leq B$, where $B$ is fixed to be $1$.
We also sample $f_{1}$ and $f_{2}$ such that both $f_{1}(0) > 0$ and $f_{2}(0) > 0$ so that there is at least one point in the safe set for both functions.
In our experiments, we consider one change point at $t_{c} = 150$. 
The time horizon is assumed to be $300$.
In the experiment shown below, we draw $500$ samples of function pairs $f_{1}$ and $f_{2}$.
For each pair of functions, the initial safe seed is the same for Adaptive-SafeOpt and Genie-CP; also the safe-seed set is the same for Genie-CP-SS and Genie-SS.
In Figure \ref{fig:pseudo regert with outlier} we illustrate $\Delta_{\pi}(t)$ for the different algorithms as a function of time.
We plot the average of $\Delta_{\pi}(t)$ over the $500$ samples of $(f_{1}, f_{2})$ with the standard deviation around the mean.
We observe that GP-UCB-CP, Genie-CP-SS, and Genie-SS converge to the minimum possible $\Delta_{\pi}(t)$ after an initial exploration phase.
\begin{figure}[ht!]
    \centering
    \includegraphics[width=90mm, height=45mm]{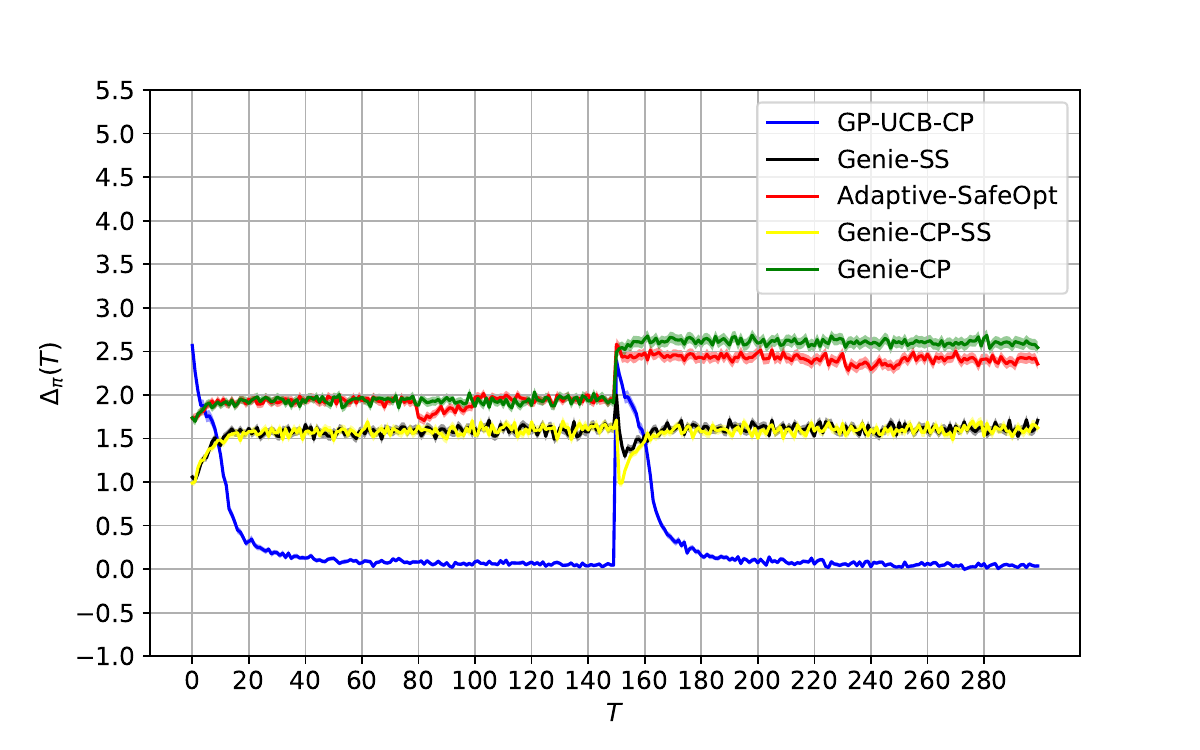}
    \caption{Comparison of $\Delta_{\pi}(T)$ as a function of $T$ for different algorithms. The change point $t_{c} = 150$. Observations are noisy with noise variance of $0.2$.}
    \label{fig:pseudo regret with outlier and noise}
\end{figure}
We also note that the proposed Adaptive-SafeOpt as well as Genie-CP converges but since the safe set that they explore is limited in size, the convergence is to a local maxima. 
The FixedWindow-SafeOpt algorithm is observed not to converge.
We also consider a case with observation noise variance of $0.2$ in Figure \ref{fig:pseudo regret with outlier and noise}.
We observe that in this case Genie-CP-SS and Genie-SS are limited by their ability to explore the safe sets completely and have larger gaps from the optimal value, in comparison to the GP-UCB-CP algorithm which is able to achieve $\Delta_{\pi}(T) = 0$ on average.
Again, the proposed Adaptive-SafeOpt converges to the local maxima corresponding to the safe seed that it finds, which is shown by the match with the Genie-CP policy.

\begin{figure}[ht!]
    \centering
    \includegraphics[width=90mm, height=45mm]{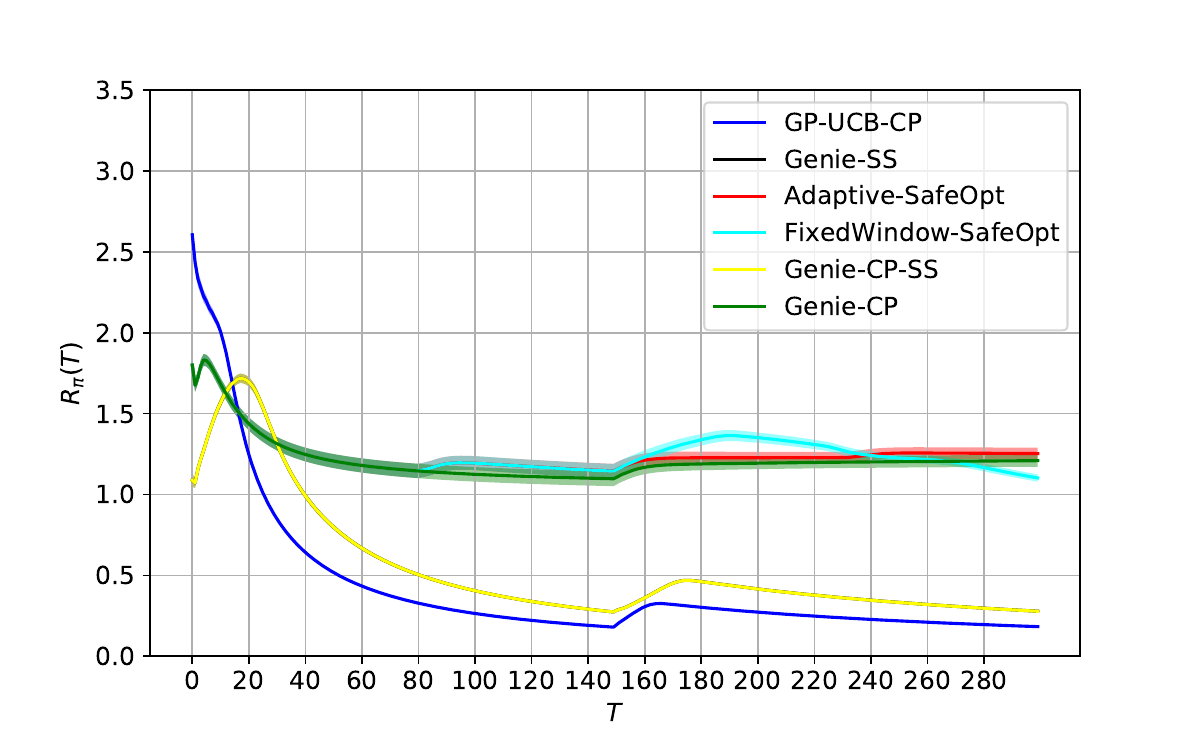}
    \caption{Comparison of $R_{\pi}(T)$ as a function of $T$ for the different algorithms. The change point $t_{c} = 150$.} 
    \label{fig:norm regert with outlier}
\end{figure}

\begin{figure}
    \centering
    \includegraphics[width=90mm, height=45mm]{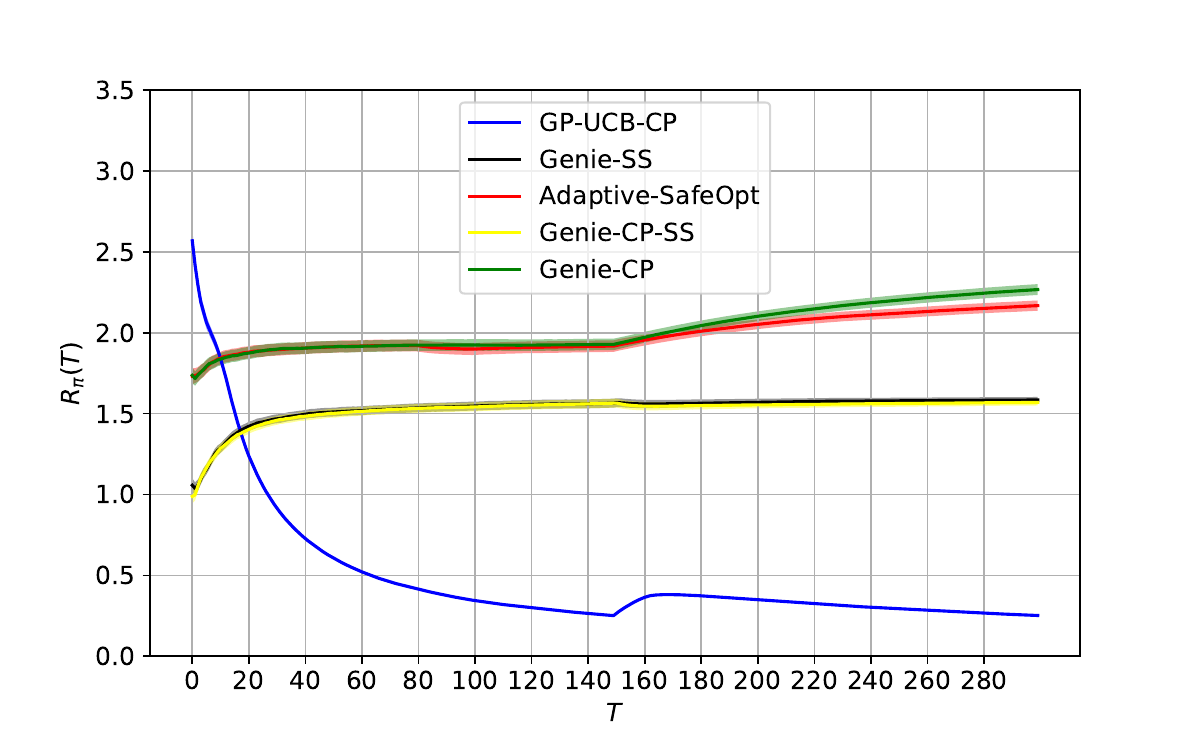}
    \caption{Comparison of $R_{\pi}(T)$ as a function of $T$ for the different algorithms. The change point $t_{c} = 150$. Observations are noisy with noise variance of $0.2$.} 
    \label{fig: norm regret with obs var}
\end{figure}

\begin{figure}[ht!]
    \centering
    \includegraphics[width=90mm, height=45mm]{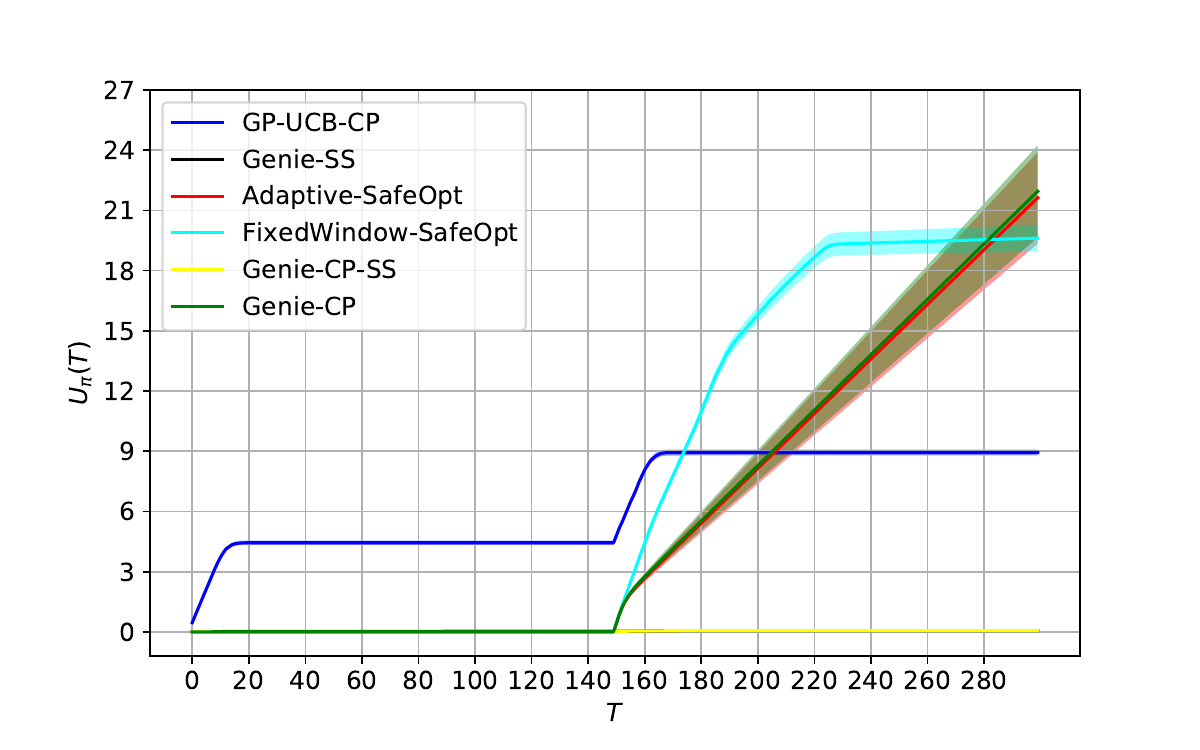}
    \caption{Comparison of the cumulative number of unsafe evaluations $U_{\pi}(T)$ as a function of $T$ for different algorithms. The number of unsafe evaluations increase at the change point $t_{c}$.}
    \label{fig:unsafe eval with outlier}
\end{figure}

\begin{figure}[ht!]
    \centering
    \includegraphics[width=90mm, height=45mm]{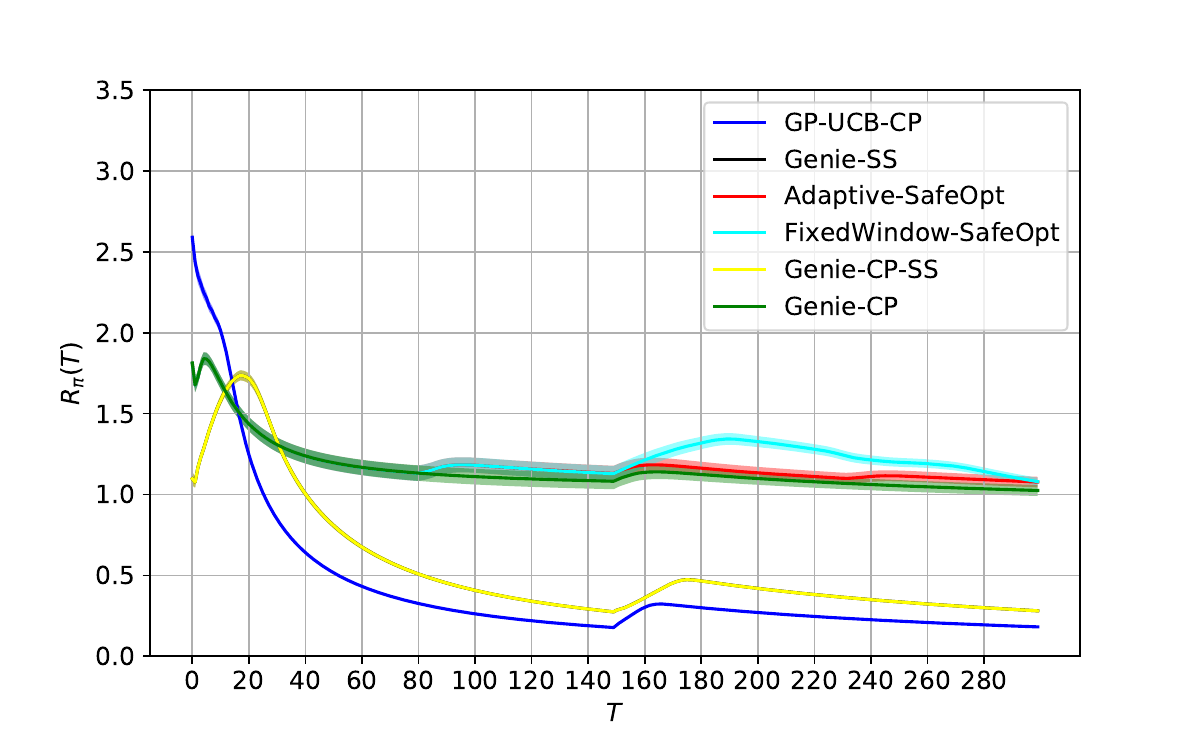}
    \caption{Comparison of $R_{\pi}(T)$ as a function of $T$ for the different algorithms after excluding the cases in which the local maxima occurs. The change point $t_{c} = 150$.} 
    \label{fig:norm regret without outlier}
\end{figure}

\begin{figure}[ht!]
    \centering
    \includegraphics[width=90mm, height=45mm]{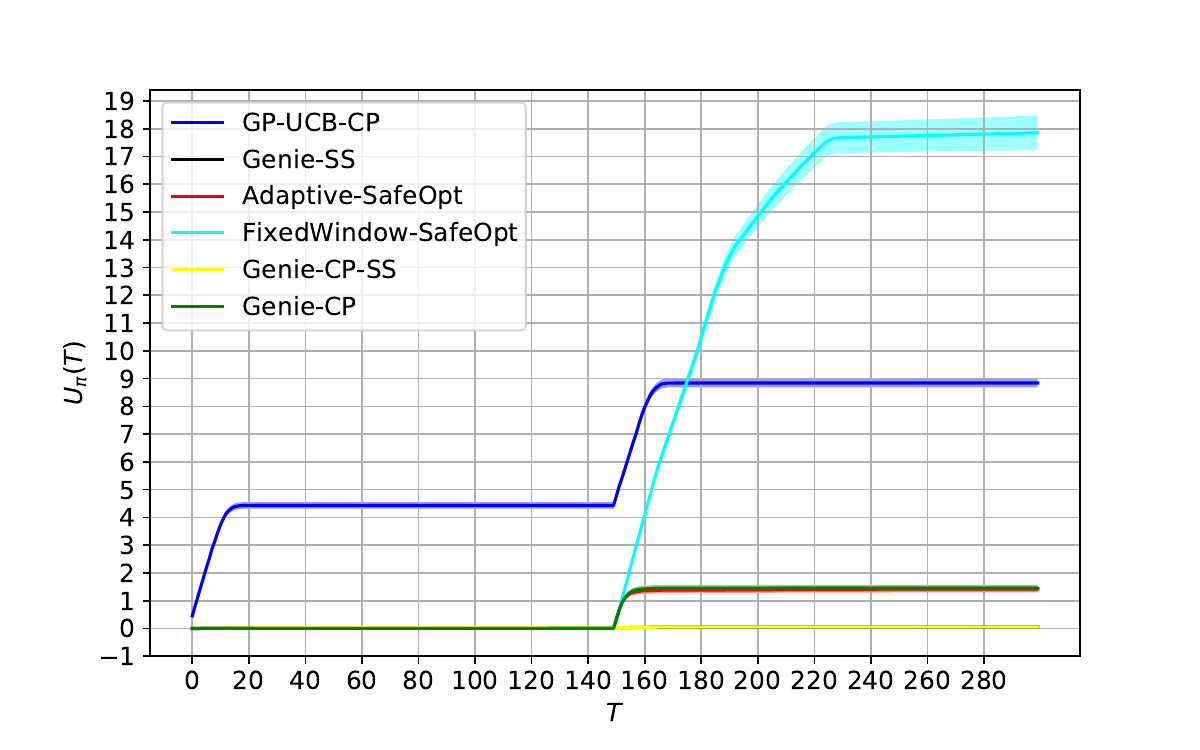}
    \caption{Comparison of the cumulative number of unsafe evaluations $U_{\pi}(T)$ as a function of $T$ for different algorithms. The cases where the local maxima occurs are excluded in the average.}
    \label{fig:unsafe eval without outlier}
    \vspace{-0.2in}
\end{figure}

The time normalized regret $R_{\pi}(T)$ for these policies without and with observation noise variance are shown in Figures \ref{fig:norm regert with outlier} and \ref{fig: norm regret with obs var}.
We observe that Genie-CP-SS and Genie-SS overlap with each other due to zero observation noise variance.
Also, in this case, Genie-SS able to detect the change point accurately at $t_{c}$ without any delay.
We illustrate the cumulative number of unsafe evaluations $U_{\pi}(T)$ for the different policies in Figure \ref{fig:unsafe eval with outlier}.
Interestingly, we find that on average, $U_{\pi}(T)$ increases for those policies for which the side information about the safe set is not available.
This is found to happen because the proposed algorithms get attracted to local maxima, which are unsafe.
Another set of experiments where the averaging is done by excluding such examples confirm this; see Figures \ref{fig:norm regret without outlier} and \ref{fig:unsafe eval without outlier}.
We then observe that GP-UCB has traded off unsafe evaluations with achieving the global maxima.
We note that the performance of Adaptive-SafeOpt depends critically on the safe-set initialization at the change point.
Although Adaptive-SafeOpt is able to converge to a local safe maxima, it could still be larger than the global safe maxima which is achieved by Genie-CP-SS or Genie-SS.
It has also been observed that instead of choosing $x_{t + 1}$ as $\argmax l_{t}(x)$ when $S_{t} = \emptyset$, the GP-UCB choice of $x_{t + 1} = \argmax u_{t}(x)$ leads to lower $R_{\pi}(T)$.
\vspace{-2em}
\section{Conclusion and Future Work}
In this paper, we considered the problem of safe optimization in a switching environment using the framework of Bayesian optimization and change point detection.
We proposed a heuristic algorithm called Adaptive-SafeOpt for this purpose and evaluated the performance of the algorithm via simulations.
We observed that a major challenge in extending safe optimization to switching environments is finding a safe point to continue exploration at the change time.
In future, we plan to extend this to the MDP setting and also to obtain worst case as well as instance specific lower bounds to the performance of such safe optimization algorithms.
An extensive study of the performance of the proposed algorithms as a function of the parameters is also planned.
\bibliographystyle{ieeetr}
\bibliography{ref}
\end{document}